\documentclass[11pt]{article}
\newcommand{\proof}{\noindent {\bf Proof: }}
\newcommand{\remark}{\noindent {\bf Remark: }}
\newcommand{\corollary}{\noindent {\bf Corollary: }}

\newtheorem{theorem}{Theorem}
\newtheorem{lemma}{Lemma}

\def\qed{\hfill $\Box$}

\usepackage{amssymb}
\usepackage{graphicx}

\begin{document}
\title{Formulas on hyperbolic volume\footnote{Dedicated to the memory of J.Bolyai}}
\author{\'A.G.Horv\'ath\\ Department of Geometry, \\
Budapest University of Technology and Economics,\\
H-1521 Budapest,\\
Hungary\\
e-mail: ghorvath@math.bme.hu}
\date{October 15, 2010}

\maketitle

\begin{abstract}
This paper collects some important formulas on hyperbolic volume. To determine concrete values of volume function is a very hard question requiring the knowledge of various methods. Our goal to give a non-elementary integral on the volume of the orthosceme (obtain it without using the Schl\"afli differential formula), using edge-lengthes as the only parameters.
\end{abstract}

{\bf MSC(2000):} 51F10, 52B10

{\bf Keywords:} coordinate systems, formulas on hyperbolic volume, Lobachevsky integral, orthosceme

\section{Introduction}

In the first section we give certain formulas with respect to some important coordinate systems and models, respectively.

Then we collect the classical results on three-dimensional hyperbolic volume of J.Bolyai and N.I.Lobachevsky. The most famous volume-integral (dependent on the dihedral angles of the orthosceme) discovered by N.I.Lobachevsky known and investigated worldwide however it is not well-known that for this volume J.Bolyai also gave two integrals. He used as parameters both of the measure of the dihedral angles and the edges, respectively. We observed that there is no volume-formula by edge-lengthes as parameters so as an application of our general formulas we compute such an integral. We will use to this calculation the system of hyperbolic orthogonal coordinates.

Finally we give a collection some new interesting formulas of special kinds of bodies discovered by contemporary mathematics showing that this old and hard problem is evergreen.

\subsection{Notation}

\begin{description}

\item[$\mathbb{R}^{n}$, $\mathbb{E}^{n}$, $\mathbb{H}^{n}$:] The space of the $n$-dimensional real vectors, the $n$-dimensional Euclidean space and the $n$-dimensional hyperbolic space, respectively.

\item[$x_i$:] The $i^{th}$ coordinate axis, and coordinate value with respect to an Cartesian coordinate system of  $\mathbb{E}^{n}$ or $\mathbb{H}^{n}$.

\item[$\xi_i$:] The $i^{th}$ coordinate axis, and coordinate value with respect to paracyclic coordinate system of $\mathbb{H}^{n}$.

\item[$\phi_i$, $r_i$:] The angle between the $i^{th}$ coordinate axis and the segment $OP_{i+1}$ and the length of the segment $OP_{i}$, respectively.

\end{description}

\section{General formulas}

In hyperbolic geometry, we have a good chance to get a concrete value of the volume function if we can transform our problem into either a suitable coordinate system or a model of the space, respectively. In this section, we give volume-integrals with respect to some important system of coordinates. In our computation we also use the parameter $k$ (used by J.Bolyai for the express the curvature of the hyperbolic space).

\subsection{Coordinate system
based on paracycles (horocycles)}

\begin{figure}[htbp]
\centerline{\includegraphics[scale=0.5]{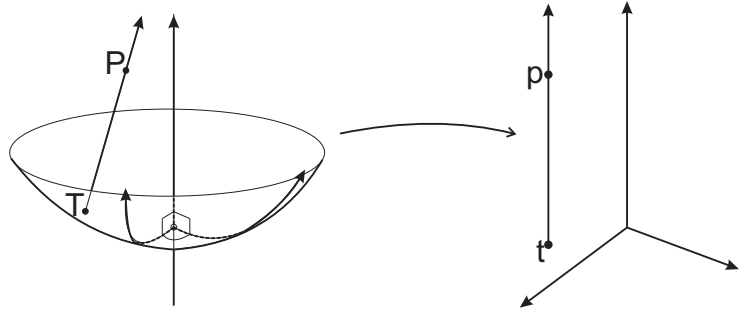}}
\caption{coordinate system based on paracycles.}
\end{figure}

Consider a parasphere of dimension $n-1$ and its pencil of rays of parallel lines. The last coordinate axis let be one of these rays, the origin will be the intersection of this line by the parasphere. The further axes are pairwise orthogonal paracycles. The coordinates of $P$ are
$
\left( \xi _1 ,\xi _2,\cdots ,\xi _n \right)^T,
$
where the last coordinate is the distance of $P$ and the parasphere, while the further coordinates are the coordinates of the orthogonal projection $T$ with respect to the Cartesian coordinate system giving by the mentioned paracycles.

In $\mathbb{R}^n$ we can correspond to $P$ a point $p$ (see on Fig.1) with ordinary Cartesian coordinates:
$$
\left( x_1,x_2,\cdots,x_n \right)^T=
\left( e^{-\frac{\xi_n}{k}}\xi_1, e^{-\frac{\xi_n}{k}}\xi_2,\cdots ,e^{-\frac{\xi_n}{k}}\xi_{n-1},\xi _n\right)^T.
$$
By definition let the volume of a Jordan measurable set $D$ be
$$
v(D):=v_n\int \limits_{D^\star}  \mathrm{d}x_1\cdots \mathrm{d}x_n,
$$
where $D^\star$ is the image of $D$ by the above mapping and $v_n$ is a constant which we will choose later. Our first formula on the volume is:
$$
v(D)=v_n\int\limits_{D} e^{-(n-1)\frac{\xi_n}{k}} \mathrm{d}\xi_1\cdots \mathrm{d}\xi_n,
$$
depending on the coordinates of the points of $D$, with the given system of coordinates.
Let now the domain $D=[0,a_1]\times \cdots\times [0,a_{n-1}]$ be a sector of parallel segments of length $a_n$ based on a coordinate-brick of the corresponding parasphere, then
$$
v(D)=v_n\int\limits_{0}^{a_1} \cdots \int\limits_{0}^{a_n} e^{-(n-1)\frac{\xi_n}{k}} \mathrm{d}\xi_n \cdots\mathrm{d}\xi_1=\frac{kv_n}{n-1}\prod \limits_{i=1}^{n-1} a_i\left[-e^{-(n-1)\frac{a_n}{k}}+e^0\right]=
$$
$$
=\frac{kv_n}{n-1}
\prod \limits_{i=1}^{n-1} a_i[1-e^{-(n-1)\frac{a_n}{k}}].
$$
If $a_n$ tends to infinity and  $a_i=1$ for $i=1\cdots (n-1)$, then the volume is equal to $\frac{kv_n}{n-1}$. Note that J.Bolyai and N.I.Lobachevski used the value $v_n=1$ so in their calculations the value of the volume is independent from the dimension but depends on the constant $k$ which determine the measure of the curvature of the space. We follows them we will determine the constant $v_n$ such that for every fixed $k$ the value of the measure of a thin layer divided by its height tends to the value of the measure of the limit figure of lower dimension. Now the limit:
$$
\lim\limits_{a_n\rightarrow \infty}\frac{v(D)}{a_n}=\frac{kv_n}{n-1}
\prod \limits_{i=1}^{n-1} a_i\lim\limits_{a_n\rightarrow \infty}\frac{[1-e^{-(n-1)\frac{a_n}{k}}]}{a_n}=v_n\prod \limits_{i=1}^{n-1} a_i,
$$
is equal to $v_{n-1}\prod \limits_{i=1}^{n-1} a_i$ showing that $1=v_1=v_2=\ldots =v_n=\ldots$.

Thus $v_n=1$ as it used by earlier. On the other hand if for a fixed $n$ the number $k$ tends to infinity the volume of a body tends to the euclidean volume of the corresponding euclidean body. In every dimension $n$ we also have a $k$ for which the corresponding hyperbolic $n$-space contains a natural body with unit volume, if $k$ is equal to $n-1$ then the volume of the paraspheric sector based on a unit cube of volume $1$ is also $1$.

So with respect to paracycle coordinate system our volume function by definition is
$$
v(D)=\int\limits_{D} e^{-(n-1)\frac{\xi_n}{k}} \mathrm{d}\xi_1\cdots \mathrm{d}\xi_n.
$$

\subsection{Volume in the Poincare half-space model}

\begin{figure}[htbp]
\centerline{\includegraphics[scale=0.5]{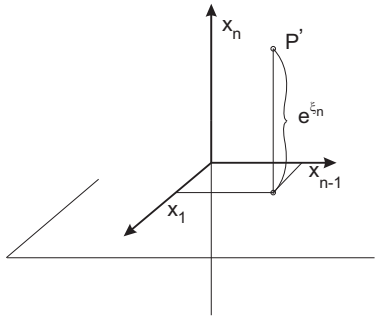}}
\caption{Coordinate system in the half-space model.}
\end{figure}

In the Poincare half-space model we consider a Cartesian coordinate system via Fig. 2. The first $n-1$ axes lie in the bounding hyperplane the last one is perpendicular to it.
If the hyperbolic coordinates  of a point $P$ (with respect to a paracycle coordinate system) is
$\left( \xi _1 ,\xi _2 ,\cdots ,\xi _n \right)^T,
$
correspond to $P$  the point $P'$ with coordinates:
$$
\left( x _1 , x _2 ,\cdots ,x _n \right)^T=
\left(\xi _1 ,\xi _2 ,\cdots ,\xi _{n-1},e^{\frac{\xi _n}{k}} \right)^T.
$$
The Jacobian of this substitution is $\frac{k}{x_n}$ showing that
$$
v(D)=k\int\limits_{D'} \frac{1}{x_n^n} \mathrm{d}x_1\cdots \mathrm{d}x_n,
$$
where $D'$ means the image $D$ via the mapping ${\bf \xi} \rightarrow {\bf x} $.

\subsection{Hyperbolic orthogonal coordinate system}

Put an orthogonal system of axes to paracycle coordinate system such that, the new half-axes let the tangent half-lines at the origin of the old one. (We can see it on Fig.3.) To determine the new coordinates of the point $P$ we project $P$ orthogonally to the hyperplane spanned by the axes $x_1,x_2,\cdots, x_{n-2},x_n$. The getting point is $P_{n-1}$. Then we project orthogonally $P_{n-1}$ onto the $(n-2)$-space spanned by $x_1,x_2,\cdots, x_{n-3},x_n$. The new point is $P_{n-2}$. Now the $(n-1)^{th}$ coordinate is the distance of $P$ and $P_{n-1}$, the $(n-2)^{th}$ is the distance of $P_{n-1}$ and $P_{n-2}$ and so on ... In the last step we get the $n^{th}$ coordinate which is the distance of the point $P_1$ from the origin $O$.
\begin{figure}[htbp]
\centerline{\includegraphics[scale=0.5]{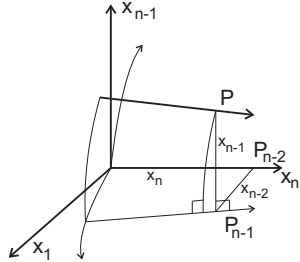}}
\caption{Coordinate system based on orthogonal axes.}
\end{figure}
Since the connection between the distance $2d$ of two points of a  paracycle and the arclength of the connecting paracycle arc $2s$ is
$$
s=k\sinh\frac{d}{k}
$$
thus the distance $z$ of the respective halving points can be calculated as:
$$
z=k\ln\cosh\frac{d}{k}.
$$
Now elementary calculation shows that the connection between the coordinates with respect to the two system of coordinates is:

\begin{eqnarray*}
\xi_{n-1} & = & e^{\frac{\xi_n}{k}}k\sinh  \frac{x_{n-1}}{k}\\
\xi_{n-2} & = & e^{\frac{\xi_n}{k}+\ln \cosh\frac{x_{n-1}}{k}}k\sinh \frac{ x_{n-2}}{k}\\
& \vdots & \\
\xi_1 & = & e^{\frac{\xi_n}{k}+\ln \cosh\frac{x_{n-1}}{k} +\cdots + \ln \cosh\frac{x_2}{k}}k\sinh \frac{ x_{1}}{k}\\
x_n & = & \xi_n+k\ln \cosh\frac{x_{n-1}}{k}+\cdots + k\ln \cosh\frac{x_2}{k} + k\ln \cosh\frac{x_{1}}{k}.
\end{eqnarray*}

From this we get a new one, corresponding point $(\xi_1,\xi_2,\cdots, \xi_n)^T\in\mathbb{H}^n$ to point  $(u_1,\cdots u_n)^T\in \mathbb{R}^n$ as in our first calculation. The corresponding system of equation is:

\begin{eqnarray*}
u_1& = &k\cosh \frac{ x_{2}}{k} \cdots  \cosh \frac{ x_{n-1}}{k}\sinh \frac { x_1}{k}\\
& \vdots & \\
u_{n-1} & = & k\sinh \frac{ x_{n-1}}{k}\\
u_n & = &  x_n-k\ln \cosh \frac{ x_{1}}{k}-\cdots - k\ln \cosh \frac{ x_{n-1}}{k}
\end{eqnarray*}
The Jacobian of this transformation is
$$
\left(\cosh \frac{ x_{1}}{k}\right)\left(\cosh  \frac{x_{2}}{k}\right)^2\cdots \left(\cosh \frac{ x_{n-1}}{k}\right)^{n-1},
$$
and we get our third formula on the volume:
$$
v(D)=\int\limits_{D} \left(\cosh \frac{ x_{n-1}}{k}\right)^{n-1}\cdots \left(\cosh  \frac{x_{2}}{k}\right)^2\left(\cosh \frac{ x_{1}}{k}\right)\mathrm{d}x_1\cdots \mathrm{d}x_n.
$$
Here we use hyperbolic orthogonal coordinates.

\subsection{Coordinate system based on spherical hyperbolic coordinates}

\begin{figure}[htbp]
\centerline{\includegraphics[scale=0.5]{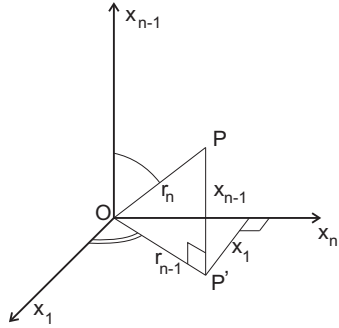}}
\caption{Spherical coordinates}
\end{figure}

From the hyperbolic orthogonal coordinates $x_1,x_2,\cdots, x_n$ we can get the spherical coordinates of a point $P$. Let the distance of $P=P_n$ from the origin is $r_n$ and denote by $\phi_i$ the angle between the $i^{th}$ coordinate axis and the segment $OP_{i+1}$ for $i=n-1,n-2,\cdots 1$. ( Here $P_{n-1}$ is the orthogonal projection of $P$ into the coordinate subspace of the axes $x_1, x_{n-2},x_n$, and so on...) We have
$$
\sinh \frac{ x_{n-1}}{k}=\sinh \frac{ r_n}{k} \cos\phi_{n-1}
$$
by the hyperbolic theorem of Sin. For general $i$, we get that
$$
\cosh \frac{ x_{n-1}}{k}\cdots \cosh \frac{x_{n-i+1}}{k}\sinh \frac {x_{n-i}}{k}=\sinh \frac{ r_n }{k}\sin\phi_{n-1}\cdots \sin\phi_{n-i+1}\cos\phi_{n-i},
$$
and by the Pythagorean theorem we can get a last equation:
$$
\cosh \frac{ x_{n-1}}{k}\cdots \cosh \frac{x_{2}}{k}\sinh \frac{  x_{1}}{k}=\sinh \frac{ r_n }{k}\sin\phi_{n-1}\cdots \sin\phi_{2}\cos\phi_{1}.
$$
Straightforward computation shows that:

$$
v(D)=k^{n-1}\int\limits_{D} \left(\sinh\frac{r_n}{k}\right)^{n-1}\sin^{n-2}\phi_{n-1}\cdots \sin\phi_{2}
\mathrm{d}\phi_1\cdots \mathrm{d}\phi_{n-1}\mathrm{d}r_n.
$$

\subsection{Volume in the projective model}

At the origin of the projective (Kayley-Clein) model we consider a Cartesian system of coordinates. Regarding the factor $k$ we assume that the radius of the sphere is $k$. The considered system is an orthogonal coordinate system both of the embedding Euclidean and the modeled hyperbolic space.
\begin{figure}[htbp]
\centerline{\includegraphics[scale=0.5]{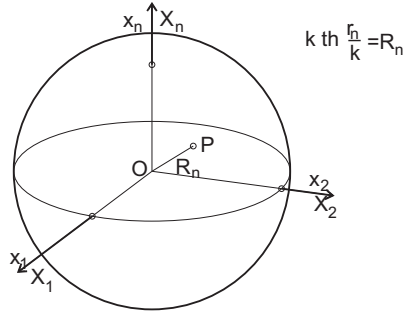}}
\caption{Coordinates in the Cayley-Klein model}
\end{figure}

We connect the hyperbolic spherical coordinates and the euclidean spherical ones by the system of equations:
$$
r_n= k\tanh^{-1}\frac{R_n}{k}, \hspace{0,5cm}
\phi_i=\theta_i, \hspace{0,5cm} i=1,\cdots, n-1.
$$
Since the Jacobian is
$$
\frac{1}{1-\left(\frac{R_n}{k}\right)^2},
$$
and
$$
\left(\sinh \frac{ r_n}{k}\right)^{n-1}=\left(\sinh\left(\frac{1}{2}\ln\frac{1+\frac{R_n}{k}}{1-\frac{R_n}{k}} \right)\right)^{n-1} =\left(\frac{\frac{R_n}{k}}{\sqrt{1-\left(\frac{R_n}{k}\right)^2}}\right)^{n-1}
$$
the volume is:
$$
v(D)=\int\limits_{D} \frac{ R_n^{n-1}}{\sqrt{1-\left(\frac{R_n}{k}\right)^2}^{n+1}}\sin^{n-2}\theta_{n-1}\cdots \sin\theta_{2}
\mathrm{d}\theta_1\cdots \mathrm{d}\theta_{n-1}\mathrm{d}R_n.
$$
Transforming it into the usual Cartesian coordinates ($X_1,\cdots ,X_n$) the new formula is:
$$
v(D)=\int\limits_{D} \frac{ 1}{\left(1-\left(\sum\limits_{i=1}^{n}\left(\frac{X_i}{k}\right)^2\right)\right)^{\frac{n+1}{2}}} \mathrm{d}X_1\cdots \mathrm{d}X_{n-1}\mathrm{d}X_n.
$$

\section{The three dimensional case}

\subsection{Formulas of J. Bolyai}

In this section $k$ is a constant giving the curvature of the hyperbolic space and the value of our constant $v_n$ is $1$. The following formulas can be found in \cite{bolyai}, \cite{bonola} and \cite{weszely}. Most of it can be easily determined using the results of the previous section. An important exception is the volume of the orthosceme, we will give a new formula for it in section 3.3.

\noindent{\bf Equidistant body:} Volume of the body determined by a disk of area p and the segments orthogonal to it with edge lengthes $t$, respectively:

$$
    v=\frac{1}{8}pk\left(e^{\frac{2q}{k}}-e^{-\frac{2q}{k}}\right) + \frac{1}{2}pq.
$$

\noindent {\bf Paraspherical sector:}
Volume of the sector of parallel half-lines intersecting orthogonally an paraspherical basic domain of area $p$:
$$
\frac{1}{2} pk.
$$

\noindent{\bf Sphere:} Volume of the sphere of radious $x$:

$$
\frac{1}{2}\pi k^3(e^{\frac{2x}{k}}-e^{\frac{-2x}{k}})-2\pi k^2x=\pi k^3\sinh\frac{2x}{k}-2\pi k^2x.
$$

\noindent{\bf Barrel:} Volume of the set of those points of $\mathbb{H}^3$, which distances from a fixed segment $AB$ of length $p$ is not greater then $q$:

$$
\frac{1}{4} \pi k^2p(e^{\frac{q}{k}}-e^{-\frac{q}{k}})^2.
$$

\noindent{\bf Orthosceme:} Volume of a special tetrahedron. Two edges $a$ and $b$ are orthogonal to each other and a third one $c$ (skew with respect to $a$) is orthogonal to the plane of $a$ and $b$.
\begin{figure}[htbp]
\centerline{\includegraphics[scale=0.5]{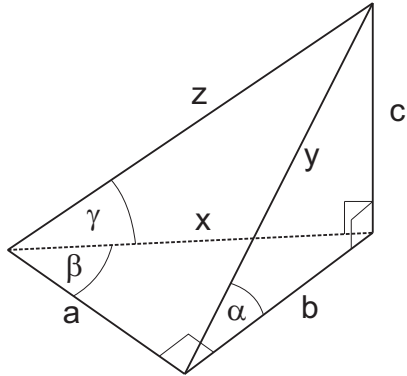}}
\caption{Orthosceme}
\end{figure}

The dihedral angle at $a$ is $\alpha$, the angle opposite to $b$ of the triangle with edges $a$ and $b$ is $\beta$ and the angle opposite to the edge $c$ in the triangle with edges $c$ and $z$ is $\gamma $, respectively (see Fig.6). J.Bolyai gave two formulas:

$$
v=\frac{\tan\gamma }{2\tan \beta}\int \limits_{0}^c\frac{z\sinh z}{\left(\frac{\cosh ^2z}{\cos ^2\alpha}-1\right)\sqrt{\frac{\cosh ^2z}{\cos ^2\gamma}-1}} {\mathrm dz},
$$

and

$$
v=\frac{1}{2}\int \limits_{0}^{\alpha}\left(-a+\frac{\sinh a\cos \phi}{2\sqrt{\tanh ^2b+\sinh ^2a\cos ^2\phi}\ln\frac{\cosh a\cos \phi +\sqrt{\tanh^2b+\sinh^2a\cos^2\phi}}{\cosh a\cos \phi -\sqrt{\tanh^2b+\sinh^2a\cos^2\phi}}}\right)
{\mathrm d\phi}
$$

\noindent{\bf Asymptotic orthosceme:} Volume of the orthosceme with ideal vertex (which is the common endpoint of the edges $a$, $x$ and $z$):

$$
v=\frac{\sin 2\alpha}{4}\int \limits_{0}^c\frac{z}{\cosh ^2z-\cos^2\alpha} {\mathrm dz}
$$

and

$$
v=\frac{1}{2}\int \limits_{0}^a\ln\frac{\cos \phi}{\sqrt{\cos ^2\phi -\tanh ^2b}} {\mathrm d\phi}
$$

\noindent{\bf Circular cone:} Volume of a cone with a basic circle of radious $b$ and with half-angle $\beta$ at its apex.

$$
v=\pi \int \limits_{0}^b\frac{\sinh ^2y}{\cosh y\sqrt{\frac{\cosh ^2y}{\cos ^2\beta}-1}} {\mathrm dy}.
$$

\noindent{\bf Asymptotic circular cone:} The apex $B$ of a circular cone tends to an ideal point on its axis of rotation.

$$
v=\pi \ln\cosh b.
$$

\subsection{Formulas of N.I.Lobachevsky}

The formulas of this subsection can be found in \cite{lobacsevszkij} or \cite{molnar}.

\noindent{\bf Barrel-wedge:} Barrel-wedge is a sector of a barrel intersected from it by two meridian-plane through its axis of rotation. Let $T$ be the area of a meridian-intersection and $p$ be the length of its parallel circular arcs. Then we have

$$
v=\frac{1}{2}pT.
$$

\noindent{\bf Orthosceme:} Let the non-rectangular dihedral angles of an orthosceme be $\alpha$, $\beta$ and $\gamma $, respectively. There admitted at the edges $a$, $z$ and $c$, respectively. (See  in Fig.6 ). Introduce the parameter $\delta$ by the equalities:
$$
\tanh \delta:=\tanh a\tan\alpha =\tanh c\tan\gamma,
$$
and the Milnor form of the Lobachevsky-function (see in \cite{milnor}):
$$
\Lambda (x)=-\int \limits_{0}^x\ln|2\sin \zeta|{\mathrm d\zeta},
$$
respectively.
Then the volume $v$ of the orthosceme is
$$
\frac{1}{4}\left[\Lambda (\alpha +\delta)-\Lambda (\alpha -\delta)-\Lambda \left(\frac{\pi}{2}-\beta +\delta\right)+\Lambda \left(\frac{\pi}{2}-\beta -\delta\right)+\Lambda (\gamma +\delta)-\Lambda (\gamma -\delta)+2\Lambda \left(\frac{\pi}{2}-\delta\right)\right].
$$

\subsection{Once more again on the volume of the orthosceme}

As an application of our general formulas we determine the volume of the orthosceme as the function of its edge-lengthes $a$, $b$ and $c$. We note that there are formulas to transform the dihedral angles into the edge-lengthes. By the notation of the previous section these are:
$$
a=\frac{1}{2}\ln\frac{\sin(\alpha +\delta)}{\sin(\alpha -\delta)}, \hspace{ 0,3cm}
c=\frac{1}{2}\ln\frac{\sin(\gamma +\delta)}{\sin(\gamma -\delta)},
\hspace{ 0,3cm}
z=\frac{1}{2}\ln\frac{\sin(\frac{\pi}{2}-\beta +\delta)}{\sin(\frac{\pi}{2}-\beta -\delta)}.
$$
It is clear that there is no simple way to get a new volume formula using these ones.

\subsubsection{The 3-dimensional case}

We now follow another way for computation, we will determine the integral
$$
v(D)=\int\limits_{D}(\cosh z)^2(\cosh y)\mathrm{d}y\mathrm{d}z\mathrm{d}y= \int\limits_{0}^{a}\int\limits_{0}^{\phi(x)}\int\limits_{0}^{\psi (x,y)}(\cosh z)^2(\cosh y)\mathrm{d}y\mathrm{d}z\mathrm{d}x,
$$
getting it from hyperbolic orthogonal coordinates using the parameter value $k=1$. The  functions $\phi (x)$ and $\psi (x,y)$ can be determined as follows. Consider the orthosceme on Fig.7.
\begin{figure}[htbp]
\centerline{\includegraphics[scale=0.5]{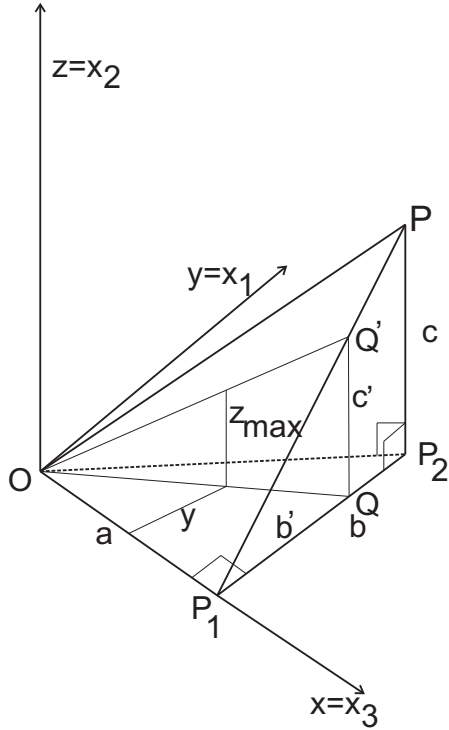}}
\caption{Orthosceme and orthogonal coordinates}
\end{figure}
In the rectangular triangle $\triangle_{OP_2P_1}$ the tangent of the angle $P_2OP_1\angle $ is:
$$
\tan P_2OP_1\angle=\frac{ \tanh b}{ \sinh a }=\frac{ \tanh y_1}{ \sinh x }.
$$
So
$$
\tanh y_{\max}=\frac{ \tanh b}{ \sinh a }\sinh x,
$$
hence
$$
0\leq y\leq \phi (x)=\tanh ^{-1}\left(\frac{ \tanh b}{ \sinh a }\sinh x\right)=:\lambda.
$$
Consider now the triangle $\triangle_{P_1P_2P_3}$. The line $O(x,y,0)$ intersects that point $Q$ for which  $|\overline{P_1Q}|=b'$, and let denote the point of the segment $PP_1$ above $Q$ be $Q'$. Thus we get the equality
$$
\tanh c'=\frac{ \tanh c}{ \sinh b }\sinh b'.
$$
Take into consideration again the equality
$$
\tanh b'=\frac{ \tanh y}{ \sinh x }\sinh a,
$$
and using the hyperbolic Pythagorean theorem, from the triangle $\triangle_{OQQ'}$ we get that
$$
\tanh z_{\max}=\tanh c'\frac{\sinh \left( \cosh ^{-1}(\cosh x\cosh y)\right)}{\sinh \left( \cosh ^{-1}(\cosh a\cosh b')\right)}=
\frac{ \tanh c}{ \sinh b }\sinh b' \frac{\sqrt{\cosh ^2x
\cosh ^2y-1}}{\sqrt{\cosh ^2a
\cosh ^2b'-1}}=
$$
$$
=\frac{ \tanh c}{ \sinh b }\sinh b' \frac{\sqrt{\sinh ^2y+\sinh ^2x
\cosh ^2y}}{\sqrt{\sinh ^2b'+\sinh ^2a
\cosh ^2b'}}=\frac{ \tanh c}{ \sinh b } \sinh y \frac{\sqrt{1+\sinh ^2x
\coth ^2y}}{\sqrt{1+\sinh ^2a\coth^2b'}}=
$$
$$
=\frac{ \tanh c}{ \sinh b }\sinh y.
$$
Hence the assumption
$$
0\leq z \leq \psi (x,y)= \tanh ^{-1}\left(\frac{ \tanh c}{ \sinh b }\sinh y\right)=:\nu
$$
holds if we fixed the first two variables. Thus the required volume is:
$$
v=\int\limits_{0}^{a}\int\limits_{0}^{\lambda} \int\limits_{0}^{\nu}(\cosh z)^2(\cosh y)\mathrm{d}z\mathrm{d}y\mathrm{d}x=
$$
$$
=\int\limits_{0}^{a}\int\limits_{0}^{\lambda} \frac{1}{2}\left[z+\frac{1}{2}(\sinh 2z)\right]_0^\nu(\cosh y)\mathrm{d}y\mathrm{d}x.
$$
Using now the equality
$$
\nu =\frac{1}{2}\ln \frac{\sinh b +\tanh c\sinh y}{\sinh b -\tanh c\sinh y},
$$
we get that
$$
v=\frac{1}{4}\left\{\int\limits_{0}^{a}\left( \int\limits_{0}^{\lambda}\ln \frac{\sinh b +\tanh c\sinh y}{\sinh b -\tanh c\sinh y}\cosh y \mathrm{d}y+\right.\right.
$$
$$
\left.\left.+\int\limits_{0}^{\lambda}\sinh\left(\ln \frac{\sinh b +\tanh c\sinh y}{\sinh b -\tanh c\sinh y}\right)\cosh y \mathrm{d}y\right)\mathrm{d}x\right\}.
$$
To determine the second integral we use that $\sinh x=\frac{e^x-e^{-x}}{2}$. Now
$$
\int\limits_{0}^{\lambda}\sinh\left(\ln \frac{\sinh b +\tanh c\sinh y}{\sinh b -\tanh c\sinh y}\right)\cosh y \mathrm{d}y=
$$
$$
=\frac{1}{2}\int\limits_{0}^{\lambda} \left(\frac{\sinh b +\tanh c\sinh y}{\sinh b -\tanh c\sinh y}-\frac{\sinh b -\tanh c\sinh y}{\sinh b +\tanh c\sinh y}\right)\cosh y \mathrm{d}y=
$$
$$
=2\int\limits_{0}^{\lambda}\frac{\sinh y\cosh y}{\frac{\sinh b}{\tanh c}-\frac{\tanh c}{\sinh b}\sinh ^2 y}\mathrm{d}y
=2\int\limits_{0}^{\lambda}\frac{\sinh 2y}{2\frac{\sinh b}{\tanh c}-\frac{\tanh c}{\sinh b}\cosh 2 y+\frac{\tanh c}{\sinh b} }\mathrm{d}y=
$$
$$
=-\frac{\sinh b}{\tanh c}\left[\ln\left(2\frac{\sinh b}{\tanh c}-\frac{\tanh c}{\sinh b}\cosh 2 y+\frac{\tanh c}{\sinh b}\right)\right]_0^\lambda=
$$
$$
=-\frac{\sinh b}{\tanh c}\ln\left(2\frac{\sinh b}{\tanh c}-\frac{\tanh c}{\sinh b}\cosh 2 \lambda+\frac{\tanh c}{\sinh b}\right)+\frac{\sinh b}{\tanh c}\ln\left(2\frac{\sinh b}{\tanh c}\right).
$$
From the definition of $\lambda$ we can calculate $\cosh 2\lambda$ and get:
$$
\cosh 2 \lambda =\frac{1}{2}\left(\frac{\sinh a+\tanh b \sinh x} {\sinh a-\tanh b \sinh x}+\frac{\sinh a-\tanh b \sinh x} {\sinh a+\tanh b \sinh x}\right),
$$
thus the value of the second integral (denoted by $\mathrm{ II }$) is:
$$
\mathrm{ II }:=-\frac{\sinh b}{\tanh c}\ln\left(1-\frac{\tanh ^2c\sinh ^2x}{\cosh ^2b(\sinh ^2a-\tanh ^2b\sinh ^2x)}\right).
$$
The first part can be integrated as follows:
$$
\int\limits_{0}^{\lambda}\ln \frac{\sinh b +\tanh c\sinh y}{\sinh b -\tanh c\sinh y} \cosh y \mathrm{d}y=\left\{\left[\ln \frac{\sinh b +\tanh c\sinh y}{\sinh b -\tanh c\sinh y}\sinh y\right]_0^\lambda-\right.
$$
$$
\left.-\int\limits_{0}^{\lambda} \frac{\tanh c\cosh y[(\sinh b -\tanh c\sinh y)+(\sinh b +\tanh c\sinh y)]}{\sinh^2 b -\tanh ^2c\sinh ^2y}\sinh y \mathrm{d}y\right\}=
$$
$$
=\left\{\sinh \lambda\ln \frac{\sinh b +\tanh c\sinh \lambda}{\sinh b -\tanh c\sinh \lambda}-\int\limits_{0}^{\lambda} \frac{2\tanh c\sinh b\cosh y\sinh y}{\sinh^2 b -\tanh ^2c\cosh ^2y+\tanh ^2c}\mathrm{d}y\right\}=
$$
$$
=\left\{\sinh \lambda\ln \frac{\sinh b +\tanh c\sinh \lambda}{\sinh b -\tanh c\sinh \lambda}+\frac{\sinh b}{\tanh c}\left[\ln(\sinh^2 b -\tanh ^2c\cosh ^2y+\tanh ^2c)\right]_0^\lambda\right\}=
$$
$$
=\left\{\sinh \lambda\ln \frac{\sinh b +\tanh c\sinh \lambda}{\sinh b -\tanh c\sinh \lambda}+\frac{\sinh b}{\tanh c}\left(\ln(\sinh^2 b -\tanh ^2c\sinh ^2\lambda)-\ln(\sinh ^2b)\right)\right\}.
$$
Since
$$
\sinh ^2\lambda =\frac{\tanh ^2b\sinh ^2x}{\sinh^2 a -\tanh ^2b\sinh ^2x},
$$
the first integral is:
$$
\left\{\sinh \lambda\ln \frac{\sinh b +\tanh c\sinh \lambda}{\sinh b -\tanh c\sinh \lambda}+\frac{\sinh b}{\tanh c}\ln\left(1 -\frac{\tanh ^2c\sinh ^2x}{\cosh ^2b(\sinh^2 a -\tanh ^2b\sinh ^2x)}\right)\right\}=
$$
$$
=\sinh \lambda\ln \frac{\sinh b +\tanh c\sinh \lambda}{\sinh b -\tanh c\sinh \lambda}-\mathrm{ II }.
$$
The sum of the two parts is:
$$
\sinh \lambda\ln \frac{\sinh b +\tanh c\sinh \lambda}{\sinh b -\tanh c\sinh \lambda}.
$$

From $\lambda =\tanh ^{-1}\left(\frac{\tanh b}{\sinh a}\sinh x\right)$ follows
$$
x=\sinh^{-1}\left(\frac{\tanh\lambda\sinh a}{\sinh b}\right)=\ln\frac{\tanh\lambda\sinh a+\sqrt{\tanh^2\lambda\sinh^2a+\tanh^2b}}{\tanh b}
$$
and we get that
$$
v=\frac{1}{4}\int\limits_{0}^{b}\frac{\tanh \lambda \sinh a}{\sqrt{\tanh^2b\cosh^2\lambda+\sinh^2a\sinh^2\lambda}}\ln\left(\frac{\sinh b+\tanh c\sinh\lambda}{\sinh b-\tanh c\sinh\lambda}\right)\mathrm{d}\lambda,
$$
proving the theorem:

\begin{theorem}
Let the edges of an orthosceme be $a,b,c$, respectively where $a\bot b$ and $(a,b)\bot c$. If $k=1$ then the volume of it is:
$$
v=\frac{1}{4}\int\limits_{0}^{b}\frac{\tanh \lambda \sinh a}{\sqrt{\tanh^2b\cosh^2\lambda+\sinh^2a\sinh^2\lambda}}\ln\left(\frac{\sinh b+\tanh c\sinh\lambda}{\sinh b-\tanh c\sinh\lambda}\right)\mathrm{d}\lambda.
$$
\end{theorem}

\corollary
This formula can be simplified in the case of asymptotic orthoscemes. If the edge-length $a$ tends to infinity, the function
$
\frac{\tanh \lambda \sinh a}{\sqrt{\tanh^2b\cosh^2\lambda+\sinh^2a\sinh^2\lambda}}
$
tends to $\frac{1}{\cosh\lambda}$ showing that the volume of the orthosceme with one ideal vertex is
$$
v=\frac{1}{4}\int\limits_{0}^{b}\frac{1}{\cosh\lambda}\ln\left(\frac{\sinh b+\tanh c\sinh\lambda}{\sinh b-\tanh c\sinh\lambda}\right)\mathrm{d}\lambda.
$$
If now the length of the edge $c$ also grows to infinity, then this formula simplified into:
$$
v=\frac{1}{4}\int\limits_{0}^{b}\frac{1}{\cosh\lambda}\ln\left(\frac{\sinh b+\sinh\lambda}{\sinh b-\sinh\lambda}\right)\mathrm{d}\lambda,
$$
which is the volume of an orthosceme with two ideal vertices. If now we reflect it in the face containing the edges $b$ and $c$ then we get a tetrahedron with three ideal vertices. If than we reflect the getting tetrahedron in the face containing the edges $b$ and $a$ we get another one with four ideal vertices. The volume of it is
$$
v=\int\limits_{0}^{b}\frac{1}{\cosh\lambda}\ln\left(\frac{\sinh b+\sinh\lambda}{\sinh b-\sinh\lambda}\right)\mathrm{d}\lambda.
$$
Of this tetrahedron there are two edges ($a$ and $c$) which are skew and orthogonal to each other (its common normal transversal is $b$). Since the reflection in the line of $b$ is a symmetry of this ideal tetrahedron, we can see that
there are two types of its dihedral angles, two opposite (at the edges $a$ and $c$) are equal to each other, ( say $A$); and the other four ones are also equal to each other ( say $B$). Then we have $A+2B=\pi$, and its volume by Milnor's formula is equal to
$$
v'=\Lambda (\pi-2B)+2\Lambda (B)=\Lambda (2B)+2\Lambda (B)=4\Lambda (B)+2\Lambda \left(B+\frac{\pi}{2}\right).
$$
(We used that the Lobachevsky function is odd, periodic of period $\pi$, and satisfies the identity  $\Lambda (2B)=2\Lambda (B)+2\Lambda (B+\frac{\pi}{2}).) $ We have the following connection between the two integrals:
$$
0=\int\limits_{0}^{b}\frac{1}{\cosh\lambda}\ln\left(\frac{\sinh b+\sinh\lambda}{\sinh b-\sinh\lambda}\right)\mathrm{d}\lambda+2\int \limits_{0}^{B+\frac{\pi}{2}}\ln|2\sin \zeta|{\mathrm d\zeta}+4\int \limits_{0}^B\ln|2\sin \zeta|{\mathrm d\zeta}.
$$

\remark If we substitute to this formula the first-order terms of the Taylor series of the functions in the integrand, respectively, we get that
$$
v=\frac{1}{4}\int\limits_{0}^{b}\frac{\tanh \lambda \sinh a}{\sqrt{\tanh^2b\cosh^2\lambda+\sinh^2a\sinh^2\lambda}} \ln\left(\frac{\sinh b+\tanh c\sinh\lambda}{\sinh b-\tanh c\sinh\lambda}\right)\mathrm{d}\lambda=
$$
$$
=\frac{1}{2}\int\limits_{0}^{b}\frac{\lambda a}{\sqrt{b^2+a^2\lambda^2}} \frac{c\lambda}{b}\mathrm{d}\lambda =\frac{ac}{2b^2}\int\limits_{0}^{b}\frac{\lambda ^2 }{\sqrt{1}}\mathrm{d}\lambda=\frac{abc}{6}.
$$
This shows that our formula for little values gives back the euclidean one.

\subsubsection{The case of dimension two}

The following calculation shows that area of a rectangular triangle also can be get with our method.
$$
\int\limits_{0}^{a}\int\limits_{0}^{\phi(x)}(\cosh y)\mathrm{d}y\mathrm{d}x=\int\limits_{0}^{a}\frac{\frac{\tanh b}{\sinh a}\sinh x}{\sqrt{1-\left(\frac{\tanh b}{\sinh a}\sinh x\right)^2}}\mathrm{d}x=
$$
$$
=\frac{1}{2}\frac{\sinh a}{\tanh b}\int\limits_{0}^{(\tanh b)^2}\frac{1}{\sqrt{1+\left(\left(\frac{\sinh a}{\tanh b}\right)^2-1\right)t-\left(\frac{\sinh a}{\tanh b}\right)^2t^2}}\mathrm{d}t
$$

At this point we can use the following result on antiderivative: if $y=ax^2+bx+c$ and $a<0$ then $\int\frac{1}{\sqrt{y}}\mathrm{d}t=\frac{1}{\sqrt{-a}} \sin^{-1}\frac{-y'}{\sqrt{b^2-4ac}}$. Thus the required area is:
$$
\frac{1}{2}\left[\sin^{-1}\frac{(2\sinh^2a)t-\sinh^2a+\tanh ^2b}{\sinh^2a+\tanh ^2b}\right]^{\tanh ^2b}_0=
\frac{1}{2}\left[\sin^{-1}(2\sin^2\beta\sinh^2a-\cos2\beta)- \sin^{-1}(-\cos2\beta)\right].
$$
Using now the equalities $\sin \beta=\frac{\sinh b}{\sinh c}$, $\sin\alpha=\frac{\sinh a}{\sinh c}$ we get that $\frac{-\sinh^2a+\tanh ^2b}{\sinh^2a+\tanh ^2b}=-\cos(2\beta)$, $\frac{(2\sinh^2a)\tanh ^2b}{\sinh^2a+\tanh ^2b}=2\sin ^2\beta\sinh^2a=2\frac{\sinh ^2b}{\sinh ^2c}\sinh^2a$ and $1=\frac{\sinh ^2b\sinh^2a}{\sinh ^2c}+\frac{\sinh ^2b+\sinh^2a}{\sinh ^2c}$. Now our formula simplified into the form
$$
\frac{1}{2}\left[\sin^{-1}(2-2\frac{\sinh ^2b+\sinh^2a}{\sinh ^2c}+\sin^2\beta-\cos^2\beta)+(\frac{\pi}{2}-2\beta)\right]=\pi-\left(\alpha +\beta +\frac{\pi}{2}\right),
$$
as we stated.

\subsubsection{The case of dimension n}

The following lemma plays an important role in the $n$-dimensional case.

\begin{lemma}
We have two $k$-dimensional affine subspaces $H_k$ and $H'_k$, respectively for which they intersection have dimension $k-1$. Assume that the points $P\in H_k$, $P'\in H'_k$ and $P''\in H_k\cap H'_k$ hold the relations $PP'\bot H'_k$ and $P'P''\bot H_k\cap H'_k$, respectively. Then the angle
$$
\alpha =\tan ^{-1}\frac{\tanh(PP')}{\sinh P'P''},
$$
is independent from the position of $P$ in $H_k$.
\end{lemma}

\proof
Let $P$ and $Q$ be two arbitrary points of $H_k$. Then by the theorem of three perpendiculars we also have that $PP''\bot H_k\cap H'_k$ and $Q'Q''\bot H_k\cap H'_k$ implying that $P''Q''$ orthogonal to the lines $PP',P'P'',QQ'$ and $Q'Q''$, respectively. Thus the angles $PP'P''\angle$, $QQ'Q''\angle$ are equals. This immediately implies the statement.
\qed

Let now our orthosceme is the convex hull of the vertices $0=A_0, A_1, \ldots , A_{n}$ situated into a hyperbolic orthogonal coordinate system as we did it in the three dimensional case. More precisely, the coordinates of the vertices are $(0,\ldots, 0)^T$, $(0,\ldots, 0,a_n)^T$, $(a_1,0,\ldots, 0,a_n)^T$, $(a_1,a_2,0\ldots, 0,a_n)^T$ $\ldots $ $(a_1,a_2,\ldots,a_n)^T$, respectively. Introduce the function giving the upper boundary of the successive integrals. These are $\phi_0=a_n$, $\phi_1:x_n\mapsto x_1$ for a point  $(x_1,0,\ldots,0,x_n)^T$ of the edge $\mbox{ conv }(OA_2)$, $\phi_2:(x_n,x_1)\mapsto x_2$ on the points $(x_1,x_2,0\ldots,x_n)^T$ of the triangle $\mbox{ conv }(OA_2A_3)$.
In general $\phi_k:(x_1,\ldots,x_{k-1},0,\ldots, 0,x_n)\mapsto x_k $ if the corresponding point $(x_1,x_2\ldots x_{k-1},x_k,0,\ldots,0,x_n)^T$ is on the $k-1$-face $\mbox{ conv }(OA_2A_3\cdots A_{k+1})$ and so on...
From Lemma 1 we get that
$$
\frac{ \tanh \phi_{k+1}(x_n,x_1,\ldots ,x_k)}{\sinh x_k}=\frac{ \tanh a_{k+1}}{\sinh a_k},
$$
implying that
$$
\phi_{k+1}(x_n,x_1,\ldots ,x_k)=\tan^{-1}\left(\frac{ \tanh a_{k+1}}{\sinh a_k}\sinh x_k\right).
$$
In the case of $k=1$, the volume can be determined by the following $n$-times integral:
$
v(O)=\int\limits_{0}^{a_n}\int\limits_{0}^{\tan^{-1}\left(\frac{ \tanh a_{1}}{\sinh a_n}\sinh x_n\right)}\cdots  \int\limits_{0}^{\tan^{-1}\left(\frac{ \tanh a_{n-1}}{\sinh a_{n-2}}\sinh x_{n-1}\right)}
$

\indent \hspace{5cm} $(\cosh^{n-1}x_{n-1})(\cosh^{n-2}  x_{n-2})\cdots (\cosh  x_{1})
\mathrm{d}x_{n-1}\cdots \mathrm{d}x_{1}\mathrm{d}x_n.
$

\subsection{Further results on concrete volumes}

Here we mention some important results from the last few decades, which are calculate hyperbolic volumes. To investigate it we can see immediately that our formula on orthosceme has very simple building up. In this section the great latin letters mean the measures of the dihedral angles at the edges denoted by the corresponding small ones.

\begin{theorem}[J. Milnor \cite{milnor}]
Opposite dihedral angles of ideal tetrahedron are equal to each other, $A+B+C=\pi$ and its volume $v$ is
$$
v=\Lambda (A)+\Lambda (B)+\Lambda (C),
$$
where $\Lambda (x)=-\int \limits_{0}^x\ln|2\sin \zeta|{\mathrm d\zeta}$.
\end{theorem}

\begin{theorem}[R.Kellerhaus \cite{kellerhaus}]
Let $R$ be a Lambert cube with essential angles $w_k$, $0 \leq w_k \leq \frac{\pi}{2}, k = 0, 1, 2$.
Then the volume $V(R)$ of $R$ is given by
$$
V(R)=\frac{1}{4}\left\{\sum\limits_0^2(\Lambda(w_i+\theta)-\Lambda(w_i-\theta))- \Lambda(2\theta)+2\Lambda\left(\frac{\pi}{2}-\theta\right)\right\}
$$
with
$$
0<\theta=\tan ^{-1}\frac{\sqrt{\cosh^2V_1-\sin^2w_0\sin^2w_2}}{\cos w_0\cos w_2}\leq \frac{\pi}{2}.
$$
\end{theorem}

\begin{theorem}[Y. Mohanty \cite{mohanty}]
Let $O$ be an ideal symmetric octahedron with all vertices on the infinity. Then $C=\pi-A$, $D=\pi-B$, $F=\pi-E$ and the volume of $O$ is:
$$
v=2\left(\Lambda \left(\frac{\pi+A+B+E}{2}\right)+\Lambda \left(\frac{\pi-A-B+E}{2}\right)+\Lambda \left(\frac{\pi+A-B-E}{2}\right)+\right.
$$
$$
\left.+\Lambda \left(\frac{\pi-A+B-E}{2}\right)\right).
$$
\end{theorem}

\begin{theorem}[D. Derevin and A.Mednykh \cite{mednykh}]
The volume of the hyperbolic tetrahedron T = T (A,B,C,D,E,F) is equal to
$$
Vol(T) = \frac{-1}{4}\int\limits^{z_2}_{z_1}\log\frac{\cos\frac{ A+B+C+z}{2}\cos\frac{ A+E+F+z}{2}\cos\frac{B+D+F+z}{2}\cos\frac{C+D+E+z}{2}} {\sin\frac{A+B+D+E+z}{2}\sin\frac{A+C+D+F+z}{2}\sin\frac{ B+C+E+F+z}{2}\sin \frac{z}{2}}\mathrm{dz},
$$
where $z_1$ and $z_2$ are the roots of the integrand, given by
$$
z_1 = \tan^{-1}\frac{ k_2}{k_1}-\tan^{-1}\frac{ k_4}{k_3} \mbox{ , } z_2 = \tan^{-1}\frac{ k_2}{k_1}+\tan^{-1}\frac{ k_4}{k_3}
,
$$
with
$$
k_1 = -(\cos S + \cos (A + D) + \cos (B + E) + \cos (C + F) + \cos (D + E + F)
+
$$
$$
 \cos(D + B + C)+
+ \cos(A + E + C) + \cos(A + B + F)),
$$
$$
k_2 = \sin S + \sin (A + D) + \sin (B + E) + \sin (C + F) + \sin (D + E + F) +
$$
$$
 \sin (D + B + C)+ \sin (A + E + C) + \sin (A + B + F),
$$
$$
k_3 = 2(\sin A\sin D + \sin B \sin E + \sin C \sin F),
$$
$$
k_4 =\sqrt{k^2_1 + k_2^2- k_3^2},
$$
and $S = A + B + C + D + E + F$.
\end{theorem}

The above theorem implies the theorem of Murakami and Yano:

\begin{theorem}[J.Murakami, M.Yano \cite{murakami}]
The volume of the simplex $T$ is
$$
v=\frac{1}{2}\Im \left(U(z_1,T)-U(z_2,T)\right),
$$
where
$$
U(z,T)=\frac{1}{2}(l(z)+l(A+B+D+E+z)+l(A+C+D+F+z)+l(B+C+E++F+z)-
$$
$$
 -l(\pi+A+B+C+z)-l(\pi+A+E+F+z)-l(\pi+B+D+F+z)-l(\pi+C+D+E+z),
$$
and $l(z)=\mathrm{Li}_2(e^{iz})$ by the Dilogarithm function
$$
\mathrm{Li}_2(z)=-\int\limits_0^x\frac{\log (1-t)}{t}\mathrm{dt}.
$$
\end{theorem}

Here I mention the theorem of Y.Cho and H. Kim described the volume of a tetrahedron using Lobachevsky function, too. It is also very complicated formula the reader can find it in \cite{cho}.

\begin{center}
\'Akos G.Horv\'ath,\\
 Department of Geometry \\
Budapest University of Technology and Economics\\
1521 Budapest, Hungary
\\e-mail: ghorvath@math.bme.hu
\end{center}

\end{document}